\documentclass[final]{siamltex}

\usepackage{epsfig}
\usepackage{graphicx}
\usepackage{amsmath}
\usepackage{amssymb}

\title{Compartmental analysis of renal physiology using nuclear medicine data and statistical optimization}

\author{Sara Garbarino\thanks{Dipartimento di Matematica, Universit\`a di Genova and CNR -- SPIN, Genova, via Dodecaneso 35 16146 Genova, Italy (garbarino$@$dima.unige.it).}
\and Giacomo Caviglia\thanks{Dipartimento di Matematica, Universit\`a di Genova, via Dodecaneso 35 16146 Genova, Italy (caviglia$@$dima.unige.it).}
\and Massimo Brignone\thanks{IRCCS San Martino -- IST, Largo Benzi 10, 16132 Genova, Italy (brignone$@$dima.unige.it).}
\and Michela Massollo\thanks{IRCCS San Martino -- IST, Largo Benzi 10, 16132 Genova, Italy}
\and Gianmario Sambuceti\thanks{DISSAL, Universit\`a di Genova and IRCCS San Martino -- IST, Largo Benzi 10, 16132 Genova, Italy (sambuceti$@$unige.it).}
\and Michele Piana\thanks{Dipartimento di Matematica, Universit\`a di Genova and CRN -- SPIN, Genova, via Dodecaneso 35 16146 Genova, Italy (piana$@$dima.unige.it).}}

\begin{document}

\maketitle

\begin{abstract}
This paper describes a general approach to the compartmental modeling of nuclear data based on spectral analysis and statistical optimization. We utilize the renal physiology as test case and validate the method against both synthetic data and real measurements acquired during two micro-PET experiments with murine models.
\end{abstract}


\section{Introduction}
Positron Emission Tomography (PET) \cite{Ollinger}  is an imaging technique capable of detecting picomolar quantities of a labeled tracer with temporal resolution of the order of seconds. FDG-PET \cite{Hays,Kamasak,Qiao} is a PET modality in which [$^{18}$F]fluoro-2-deoxy-D-glucose(FDG) is used as a tracer to evaluate glucose metabolism and to detect diseases in many different organs. From a mathematical viewpoint, PET (and, specifically, FDG-PET) experiments involve two kinds of problems. First, signal processing techniques are applied to reconstruct the time dependence of location and concentration of tracer from the measured radioactivity. Second, these dynamic PET data can be processed to estimate physiological parameters that describe the functional behavior of the inspected tissues. This paper focuses on this second aspect and, specifically, examines FDG-PET data from kidneys and bladder, corrected for radiactive decay, to recover information on renal physiology. The data processed within the paper have been acquired by means of a PET device for small animals (mice), but the mathematical procedure for the treatment of data is essentially the same as for PET data of humans. 

The analysis of tracer kinetics in tissues whose activities have been measured with FDG-PET is typically based on compartmental models \cite{Schmidt-2,Gunn,Schmidt}. This conceptual framework identifies different compartments in the physiological system of interest, each one characterized by a specific (and homogeneous) functional role (therefore a compartment is not necessarily represented by a specific organ or anatomical district). As any other tracer, FDG is injected into the system with a concentration mathematically modeled by the so-called {\em{Time Activity Curve}} (TAC). The time dependent concentrations of tracer in each compartment constitute the state variables that can be determined from PET data. The time evolution of the state variables (the kinetics of the system) is modelled by a linear system of ordinary differential equations for the concentrations, expressing the conservation of tracer during flow between compartments. The (constant) micro-parameters describing the input/output rates of tracer for each compartment are called {\em{exchange coefficients}} or {\em{rate constants}}, represent the physiological parameters describing the system's metabolism and are the unknowns to be estimated. 

\begin{figure}
\begin{center}
\begin{tabular}{c}
\psfig{figure=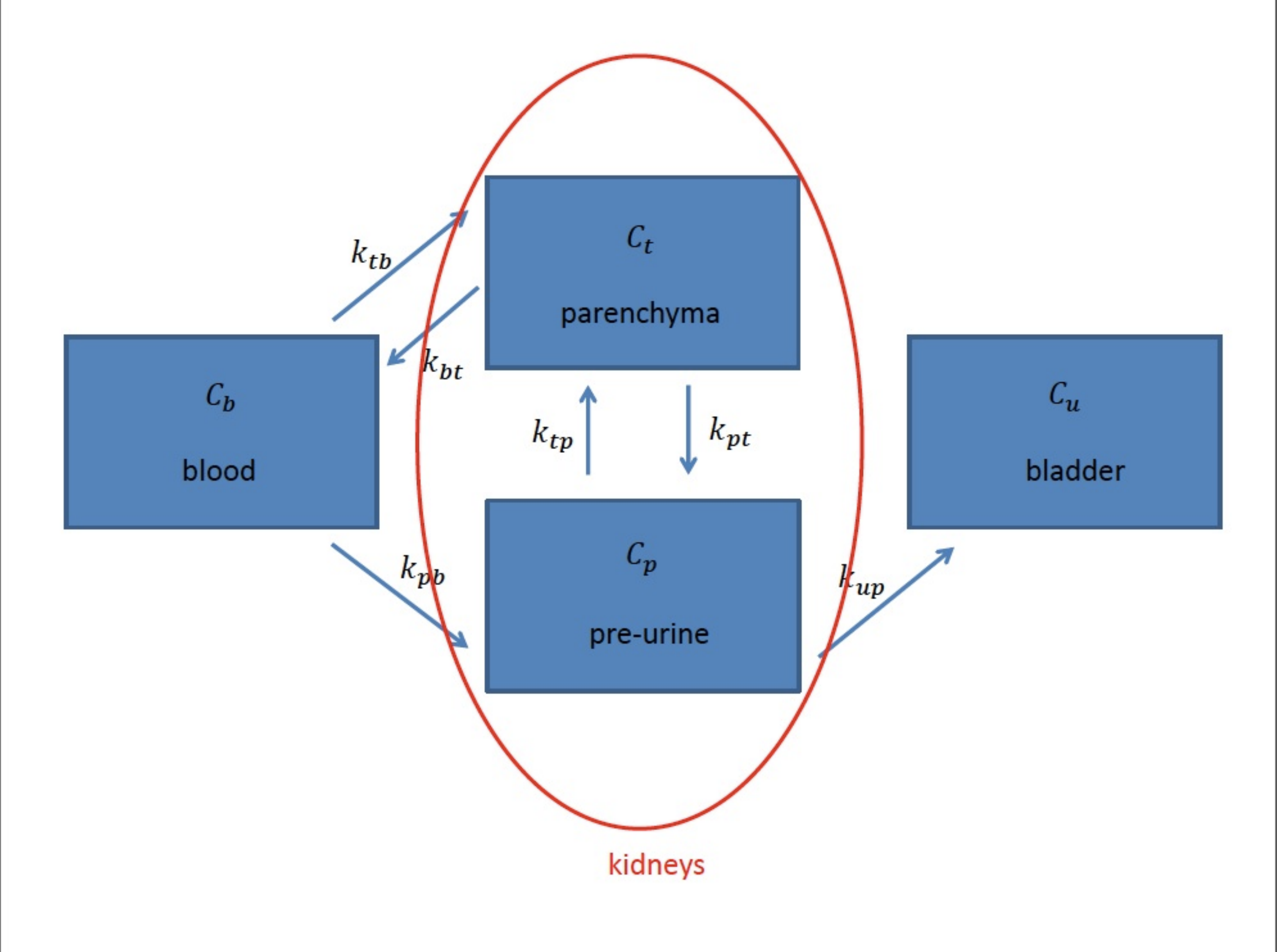,width=9.cm} 
\end{tabular}
\caption{The compartmental model of renal physiology adopted in this paper.}
\label{fig:scheme}
\end{center}
\end{figure}

The particular case of renal physiology is described in Figure \ref{fig:scheme}, where: 
\begin{itemize}
\item kidneys include two compartments, tissue (parenchyma) and pre-urine, while a third compartment, the urine, represents the bladder, where tracer is accumulated; 
\item the TAC for blood is determined from Regions of Interest (ROIs) of the left ventricle drawn on FDG-PET maps at different time steps; 
\item estimates of the tracer concentration for the bladder and for the two compartment system made of parenchyma and pre-urine are obtained by means of ROIs including bladder and kidneys, respectively; and, finally, 
\item the set of unknowns is made of six exchange coefficients, also referred to as micro parameters (other exchange coefficients, for example the ones describing the circulation for the bladder back to the kidneys, are set to zero for well-established physiological reasons).
\end{itemize}
Simpler models for tracer circulation inside kidneys can be found in the literature, involving a smaller number of exchange coefficients. While the reduction in number allows for an estimate of the unknown exchange coefficients by rather simple mathematical procedures (such as graphical methods, \cite{Logan}), this implies a notably less refined description of renal physiology.

The present paper introduces a novel method for the reduction of compartmental models like this renal one, whose main advantages are a notable generality and limited computational requirements. This method is based on spectral analysis and statistical optimization and includes two steps. In the first step, we solve the direct problem of determining the explicit formal expressions of the tracer concentration in each compartment, under the assumption that the exchange coefficients are given. According to the properties of these coefficients, solutions are classified in four different families. Consistently with a spectral analysis approach, the concentrations are expressed as linear combinations of the TAC convolved with specific exponential functions. In the second step, the inverse problem of determining the six unknown exchange coefficients is addressed by means of an 'ad hoc' implementation of Ant Colony Optimization \cite{Dorigo, Socha}, which is a statistical optimization algorithm inspired by evolutionary strategies. In this application the functional to minimize measures the discrepancy between the experimental concentrations and the analytical forms provided by the solution of the forward problem and also accounts for the blood fraction with which the vascular system of the kidneys is fed up. 

The paper content is organized as follows. In Section 2 we describe the Cauchy problem modeling the direct problem. Section 3 determines the exchange coefficients by means of statistical optimization. Section 4 shows some applications to synthetic and real PET measurements. Our conclusions are offered in Section 5.

\section{The direct problem}
The state variables of the three-compartment model adopted in this paper are the tracer concentrations in the tissue ($C_{t}$), in the pre-urine ($C_{p}$) and in the urine contained in the bladder ($C_{u}$). Moreover, the kinetic process in the system is initialized by the TAC $C_{b}$, representing the tracer concentration in blood. The usual conditions for compartmental analysis hold, e.g., tracer is uniformly distributed in each compartment at each instant \cite{Schmidt}, diffusive effects are neglected and physiological processes are in a steady state. The six constant exchange coefficients (flux rate parameters) between compartments in contact are denoted as $k_{ab}$, where the suffixes $a$ and $b$ denote the target and source compartment, respectively.  For example, $k_{up}$ is the rate constant of tracer carried ``to'' the bladder pool $u$ ``from'' the pre-urine pool $p$. We assume $k_{ab} \geq 0$ for all cases.
  
Conservation of tracer exchanged between compartments leads to the following system of linear ordinary differential equations with constant coefficients:
\begin{equation}\label{system-1}
\left\{
\begin{array}{l}
 \dot{C}_t = - (k_{bt} + k_{pt})\, C_{t} + k_{tp}\, C_{p} + k_{tb}\, C_{b} , \\
 \dot{C}_p = k_{pt} \, C_{t} - (k_{tp} + k_{up})\, C_{p} + k_{pb}\, C_{b} ,\\
  \dot{C}_u = k_{up} \, C_{p}, 
 \end{array}
 \right.
 \end{equation}
 with initial conditions $C_{t}(0)=C_{p}(0)=C_{u}(0)=0$ and where the dependence on time is implicit. The first two equations in (\ref{system-1}) can be written in the compact form
\begin{equation}\label{system-2}
\dot{\mathbf{C}} = \mathbf{A} \, \mathbf{C} + \mathbf{K} \, C_{b} , 
\end{equation}
with initial condition $\mathbf{C}(0)=0$. Here, in fact,
\begin{equation}\label{system-3}
\mathbf{C} :=  \left[ 
\begin{array}{c}
C_{t} \\ C_{p} 
\end{array}
 \right];
 \end{equation}
\begin{equation}\label{system-4}
\mathbf{A} = \left[ \begin{array}{cc} - (k_{bt}+k_{pt}) & k_{tp} \\ k_{pt} & -(k_{tp}+ k_{up}) \end{array} \right]  =: \left[ \begin{array}{cc}  -a & b \\ c & -d \end{array} \right] ;
\end{equation}
and
\begin{equation}\label{system-5}  
\mathbf{K} := \left[ \begin{array}{c} k_{tb}  \\ k_{pb}  \end{array} \right] .
\end{equation}
The definition of $\mathbf{A}$ corresponds to define 
\begin{equation}\label{system-6} 
a:= k_{bt}+k_{pt}, ~~~~~ c:= k_{pt}, ~~~~~ d:= k_{tp}+ k_{up}, ~~~~~ b:= k_{tp}
\end{equation} 
with inverse relations
\begin{equation}\label{system-7}
k_{pt} =c, ~~~~~ k_{bt}= a-c, ~~~~~ k_{tp}=b, ~~~~~ k_{up}= d-b.
\end{equation}
In the first part of this section we will show how the structure of matrix $\mathbf{A}$ is related to the explicit forms of tracer concentration in the different compartments. In the framework of compartmental analysis this issue is crucial for the study of the inverse problem performed in the next section.

Let's first assume $c \neq 0$ and $b \neq 0$. Then the solution of (\ref{system-2})-(\ref{system-7}) is given by
\begin{equation}\label{solution-full-1}
C_{t} = c_1 \,b\, E_1 + c_2 \,(d+\lambda_2) \, E_2, 
\end{equation}
\begin{equation}\label{solution-full-2}
C_{p} = c_1 \,(a+\lambda_1)\, E_1 + c_2 \,c \, E_2,
\end{equation}
\begin{equation}\label{solution-full-3}
E_i = \int_0^t e^{\lambda_i \,(t-\tau)} \, C_{b}(\tau)\,d\tau
= e^{\lambda_i \,t} * C_{b} , ~~~~~~ i=1,2, 
\end{equation}
where $*$ denotes the convolution operator, the constants $c_1$ and $c_2$ are defined as
\begin{equation}\label{solution-full-4}
c_1 = \frac {- c\,k_{tb} + (d+\lambda_2)\,k_{pb} } {(a+\lambda_1)\,(d+\lambda_2) - b\, c }, ~~~~
c_2 = \frac { (a+\lambda_1)\,k_{tb} -b \,k_{pb} } {(a+\lambda_1)\,(d+\lambda_2) - b\, c },
\end{equation}
and $\lambda_1,\lambda_2$ are the eigenvalues of the matrix ${\mathbf{A}}$ with explicit values
\begin{equation}\label{lambda-1}
\lambda_1= \frac{-(a+d) + \sqrt{(a+d)^2 - 4 (ad-bc)}}{2},
\end{equation}  
\begin{equation}\label{lambda-2}
\lambda_2= \frac{-(a+d) - \sqrt{(a+d)^2 - 4 (ad-bc)}}{2}.
\end{equation}  
We observe that both $\lambda_1$ and $\lambda_2$ are negative. The concentration $C_{u}$ is evaluated by inserting (\ref{solution-full-2}) in the third equation of (\ref{system-1}), which leads to 
\begin{equation}\label{solution-full-5}
\frac{C_{u}(t)}{k_{up}}= \int_0^t C_{p}(\sigma)\,d\sigma=\int_0^t  [ c_1\, (a +\lambda_1) \, E_1(\sigma) + c_2 \,c \,  E_2(\sigma )\big] \, d\sigma.
\end{equation}
We utilize the identity
\begin{equation}\label{identity-1}
\int_0^t  \int_0^\sigma C(\tau) \, e^{ w\,(\sigma-\tau)}d\tau d\sigma =  \frac 1 {w} \, \int_0^t C(\sigma) \, e^{ w\,(t-\sigma)} \, d\sigma  -\frac 1 w \,\int_0^t  C(\tau) \, d\tau,
\end{equation}
where $C$ is any continuous function and $w$ is a real parameter and find
\begin{equation}\label{solution-full-6}
\frac{C_{u}(t)}{k_{up}}=   
c_1\, \frac{(a +\lambda_1)}{\lambda_1}\, E_1 + c_2 \, \frac{c}{\lambda_2} \,E_2 - \big[c_1\, \frac{(a +\lambda_1)}{\lambda_1} + c_2 \, \frac{c}{\lambda_2} \big] \, \int_0^t\, C_{b}(\tau)\,d\tau.
\end{equation}

We now suppose that $b=0$ and $c \neq 0$ in the definition of the matrix $\mathbf{A}$, i.e.
\begin{equation}\label{lower-triangular}
\mathbf{A} =  \left[ \begin{array}{cc}  -a & 0 \\ c & -d \end{array} \right].
\end{equation}
In this case we set $\lambda_1 = -a$ and $\lambda_2 = -d$. This corresponds to vanishing denominators in (\ref{solution-full-4}). Integration of the linear system leads to 
\begin{equation}\label{solution-lower-triangular-1}
C_{t} = \chi_1 \,(\lambda_1-\lambda_2 ) \, E_{1}
\end{equation}
and
\begin{equation}\label{solution-lower-triangular-2}
C_{p} = c\,\chi_1 \, E_{1} + \chi_2 \, E_{2},
\end{equation}
where $E_{1}$ and $E_{2}$ are convolutions defined as in (\ref{solution-full-3}) and
\begin{equation}\label{solution-lower-triangular-3}
\chi_1 = \frac{k_{tb}}{\lambda_1-\lambda_2} , ~~~~~
\chi_2 = k_{pb}-  \frac{c \, k_{tb}}{\lambda_1-\lambda_2 }.
\end{equation}
Moreover, integrating the third equation in (\ref{system-1}) leads to
\begin{equation}\label{solution-lower-triangular-4}
\frac{C_u}{k_{up}}=   \frac{c\,\chi_1}{\lambda_1}\, E_{1,b} + \frac{\chi_2}{\lambda_2} \,E_{2,b} - \left[ \frac{c\,\chi_1}{\lambda_1} + \frac{\chi_2}{\lambda_1} \right] \, \int_0^t\, C_{b}(\tau)\,d\tau. 
\end{equation}

We then consider the case $b \neq 0$, $c=0$, i.e.
\begin{equation}\label{upper-triangular}
\mathbf{A} =  \left[ \begin{array}{cc}  -a & b \\ 0 & -d \end{array} \right]. 
\end{equation}
Again we have $\lambda_1 = -a$ and $\lambda_2 = -d$. The solution of the Cauchy problem for this upper triangular system is
\begin{equation}\label{solution-upper-triangular-1}
C_{t} = \sigma_1 \, E_{1}+ b\,\sigma_2 \, E_{2},
\end{equation}
\begin{equation}\label{solution-upper-triangular-2}
C_{p} = (\lambda_2-\lambda_1)\,\sigma_2 \, E_{2},
\end{equation}
with
\begin{equation}\label{solution-upper-triangular-3}
\sigma_1 =  k_{tb} - \frac {b\,k_{pb}}{\lambda_2-\lambda_1}  , ~~~~ \sigma_2 = \frac {k_{pb}}{\lambda_2-\lambda_1} 
\end{equation}
and $E_{1}$, $E_{2}$ defined as in (\ref{solution-full-3}).

Finally, we consider the diagonal case $b=0$, $c=0$, i.e.
\begin{equation}\label{diagonal}
\mathbf{A} = \left[ \begin{array}{cc} -a & 0 \\ 0 & -d \end{array} \right].
\end{equation}
Then, it is straightforward to obtain
\begin{equation}\label{solution-diagonal-1}
C_{t} = k_{tb} \, E_{1}, ~~~~  C_{p} = k_{pb}\, E_{2}. 
\end{equation}
and
\begin{equation}\label{solution-diagonal-2}
C_{u} = - k_{pb}\, E_{2} + k_{pb}\, \int_0^t C_{b}(\tau) \, d\tau. 
\end{equation}

\section{Solution of the inverse problem}
The model equations obtained in the previous section describe the time behavior of the tracer concentration in the three compartments of the renal system, given the TAC and the transmission coefficients. Given such equations, compartmental analysis requires the determination of the tracer coefficients by
\begin{enumerate}
\item utilizing measurements of the tracer concentrations provided by nuclear imaging, and
\item applying an optimization scheme for the solution of the inverse problem.
\end{enumerate}

In nuclear imaging experiments, the reconstructed images can provide information on the tracer concentration in the kidneys and in the bladder as well as in the input arterial blood as measured in the left ventricle. Specifically, an acquisition sequence is set up providing count data sets collected at subsequent time intervals. For each data set, an image reconstruction algorithm is applied, ROIs are drawn within the left ventricle, the kidneys and the bladder and the corresponding tracer concentrations are computed. Obviously, the tracer concentration in the kidneys is an estimate of $C_{t} + C_{p}$ plus the tracer carried by the blood contained in the kidneys' vascular system. This last term cannot be identified in the images and will be accounted for in the optimization procedure. Specifically, the optimization scheme we will implement minimizes, at each time point, the functional
\begin{equation}\label{aco-1}
{\cal{C}} := \|(C_{t} + C_{p}) - \overline{C}_{exp} \|^2 + \|C_{u} - \overline{C}_u \|^2,
\end{equation}
where $C_{u}$, $C_{t}$ and $C_{p}$ are the analytical solutions of the direct problem computed at the given time point; $\overline{C}_u$ is the concentration measured from the ROI in the bladder at the time point;  $\overline{C}_{epx} := \frac{\overline{C}_{k} - V_b C_{b}}{1-V_b}$,
where $\overline{C}_{k}$ is the concentration measured from the ROI on the kidneys, $C_{b}$ is the value of the TAC at the specific time value and $V_b$ measures the blood fraction with which the kidneys' vascular system is supplied. In the following we will assume $V_b = 0.2$, which is a physiologically sound value \cite{Bentley}. 

The minimization of functional $\cal{C}$ is realized by means of an Ant Colony Optimization (ACO) scheme \cite{Dorigo}. ACO is a statistical-based optimization method developed in the 1990s with the aim of providing a reliable although not optimal solution to some non-deterministic polynomial-time hard combinatorial optimization problems. More recently ACO has been generalized to continuos domains and has been applied to a wide range of problems. ACO takes inspiration from the way in which ants find and carry food to their nest. While an ant is going back to the nest after having taken some food, it releases a pheromone trace that serves as a trail for next ants, which are able to reach food detecting pheromone. Since the pheromone decays in time, its density is higher if the path to food is shorter and more crowded; on the other hand, more pheromone attracts more ants and at the end all ants follow the same trail. This behavior is paraphrased in ACO identifying the cost functional $\mathcal{F}$ with the length of the path to food, and the pheromone traces with a probability density which is updated at each iteration depending on the value of the cost function for a set of states. 
 In practice, at each iteration, the cost function is evaluated on a set of $P$ admissible states, and the states are ordered according to increasing values of the cost function. Then ACO defines a probability distribution which is more dense in correspondence with the cheaper states and, on its basis, $Q$ new states are extracted. A comparison procedure identifies the new best $P$ states which form the next set of states.
Formally the starting point of the algorithm is a set of $P$ states
\begin{equation}\label{aco-3}
B:=\{\mathbf{U}_k=(u_{1,k},\dots,u_{N,k})\},\ \
\end{equation}
such that
\begin{equation}\label{aco-4}
\mathbf{U}_k\in S\subset \mathbb{R}^N, \ \ \ k=1,\dots,P
\end{equation}
that are ordered in terms of growing cost, namely, $\mathcal{F}({\mathbf{U}}_1)\leq\dots\leq\mathcal{F}({\mathbf{U}}_P)$.
Next, for each $j=1,\dots,N$ and $i=1,\dots,P$, one computes the parameters
\begin{equation}\label{aco-5}
m_{i,j}=u_{j,i}, \ \ \ \ s_{i,j}=\frac{\xi}{P-1}\sum_{p=1}^P|u_{j,p}-u_{j,i}|, 
\end{equation}
and defines the probability density function
\begin{equation}\label{aco-6}
\mathcal{G}_i=\sum_{i=1}^Pw_i\mathcal{N}_{[m_{i,j},s_{i,j}]}(t),
\end{equation}
with $i=1,\dots,P$ and $\xi,q$ real positive parameters to be fixed and $w_i=\mathcal{N}_{[1,qP]}(i)$.
By sampling $S$ $Q$ times with $\mathcal{G}_j$, the procedure generates $Q$ new states $\mathbf{U}_{P+1},\dots,\mathbf{U}_{P+Q}$, enlarging the set $B$ to $\tilde{B}=\{\mathbf{U}_1,\dots,\mathbf{U}_{P+Q}\}.$ If $\mathbf{U}_{k_1},\dots,\mathbf{U}_{k_Q}$ are the $Q$ states of $\tilde{B}$ of greater cost, the updated $B$ is defined as 
\begin{equation}\label{aco-7}
B=\tilde{B}\setminus\{\mathbf{U}_{k_1},\dots,\mathbf{U}_{k_Q}\}.
\end{equation}
This procedure converges to an optimal solution of the problem by exploiting the fact that the presence of the weights $w_i$, in the definition of $\mathcal{G}_j$, gives emphasis to solutions of lower costs since $w_1>\dots>w_P.$ This fact, associated with the influence that a proper choice of parameters $\xi$ and $q$ has on the shape of the Gaussian functions, determines the way in which the method tunes the impact of the worse and best solutions.
The algorithm ends when the difference between any two states of $B$ is less than a pre-defined quantity or when the maximum allowable number of iterations is reached. The initial set $B$ of trial states is chosen by sampling a uniform probability distribution. 

The implementation of ACO for the optimization of the exchange coefficients in $\cal{C}$ is based on the following steps:
\begin{enumerate}
\item The four ACO parameters are fixed as follows. $P$ and $Q$ are chosen as in \cite{Giorgi}. Specifically, $P$ is a multiple of the number of coefficients to optimize plus one and
\begin{equation}\label{aco-8}
Q = \left[ \frac{P}{2} \right] + 1,
\end{equation}
where $[\cdot]$ indicates the floor; $q$ and $\xi$ are fixed searching for a trade off between the risk of a solution space of a too high complexity and the risk of a too high computational demand.
 \item The values of the tracer coefficients are initialized to six random numbers picked up in the interval $(0,1)$ (this choice is based on literature \cite{Shen}).
 \item The ACO procedure is then run using as $C_{u}$, $C_{t}$ and $C_{p}$ the solutions in the general full matrix case. If, during the ACO iterations, the reconstructed values for the tracer coefficients become statistically consistent with values for which the associated direct problem is described by a triangular or a diagonal matrix, the algorithm automatically switches to utilize the corresponding solution in the computation of $\cal{C}$.
\end{enumerate}
In the following section we show how this statistics-based compartmental analysis works in the case of synthetic data and real measurements recorded by a micro-PET system.

\section{Numerical experiments}
Compartmental analysis is a valid approach to physiological studies of animal models by means of PET data. An 'Albira' micro-PET system produced by Carestream Health is currently operational at the IRCCS San Martino IST, Genova and experiments with mice are currently performed by using different tracers, mainly for applications to oncology. In this section we describe the performance of our approach to compartmental analysis in the case of synthetic data simulated by mimicking 'Albira' acquisition for FDG-PET experiments. Then we will describe the results of data analysis for two real experiments, performed by using FDG.

\subsection{Application to synthetic data}
\begin{figure}
\begin{center}
\begin{tabular}{cc}
\psfig{figure=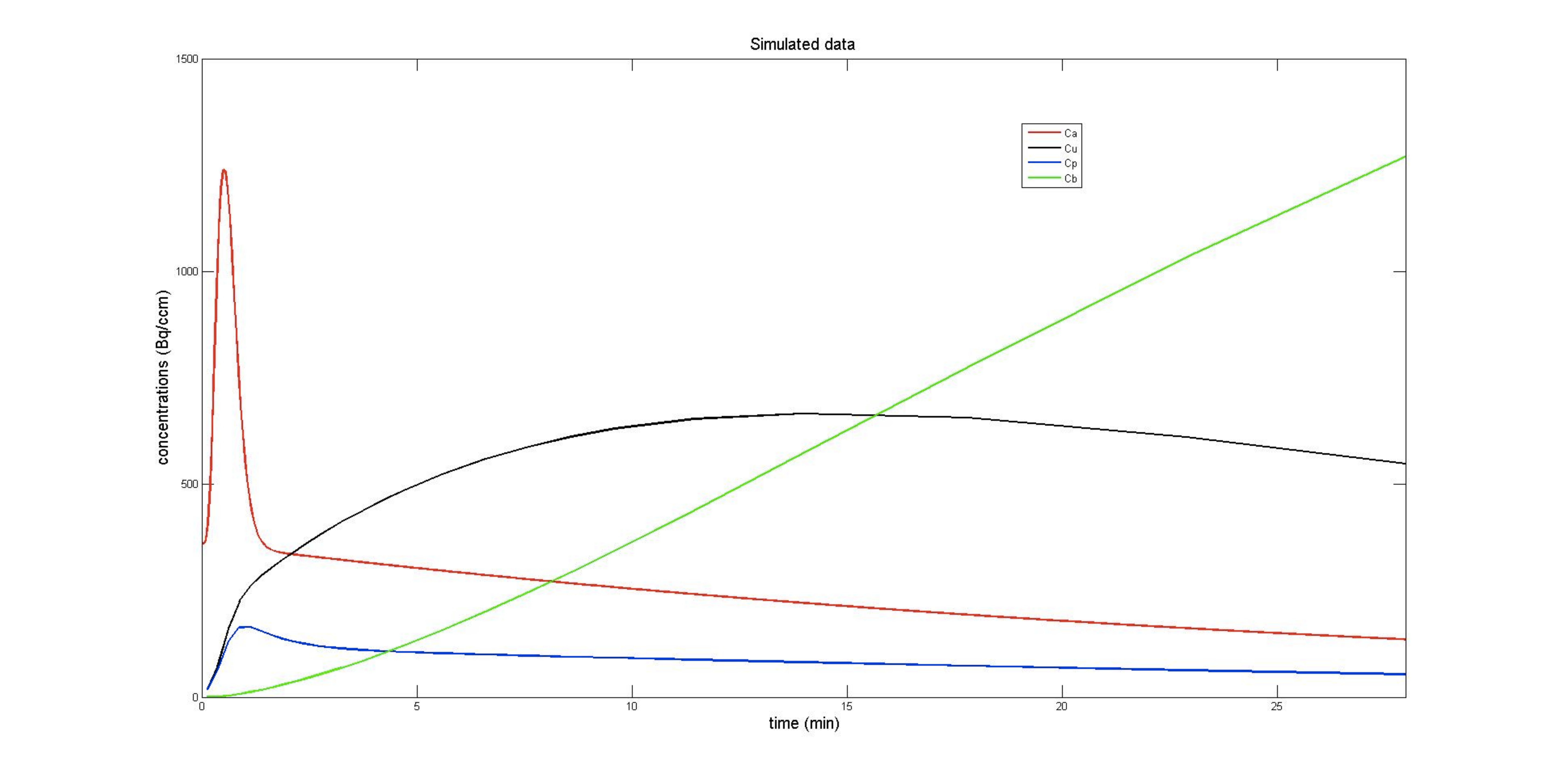,width=6.5cm} &
\psfig{figure=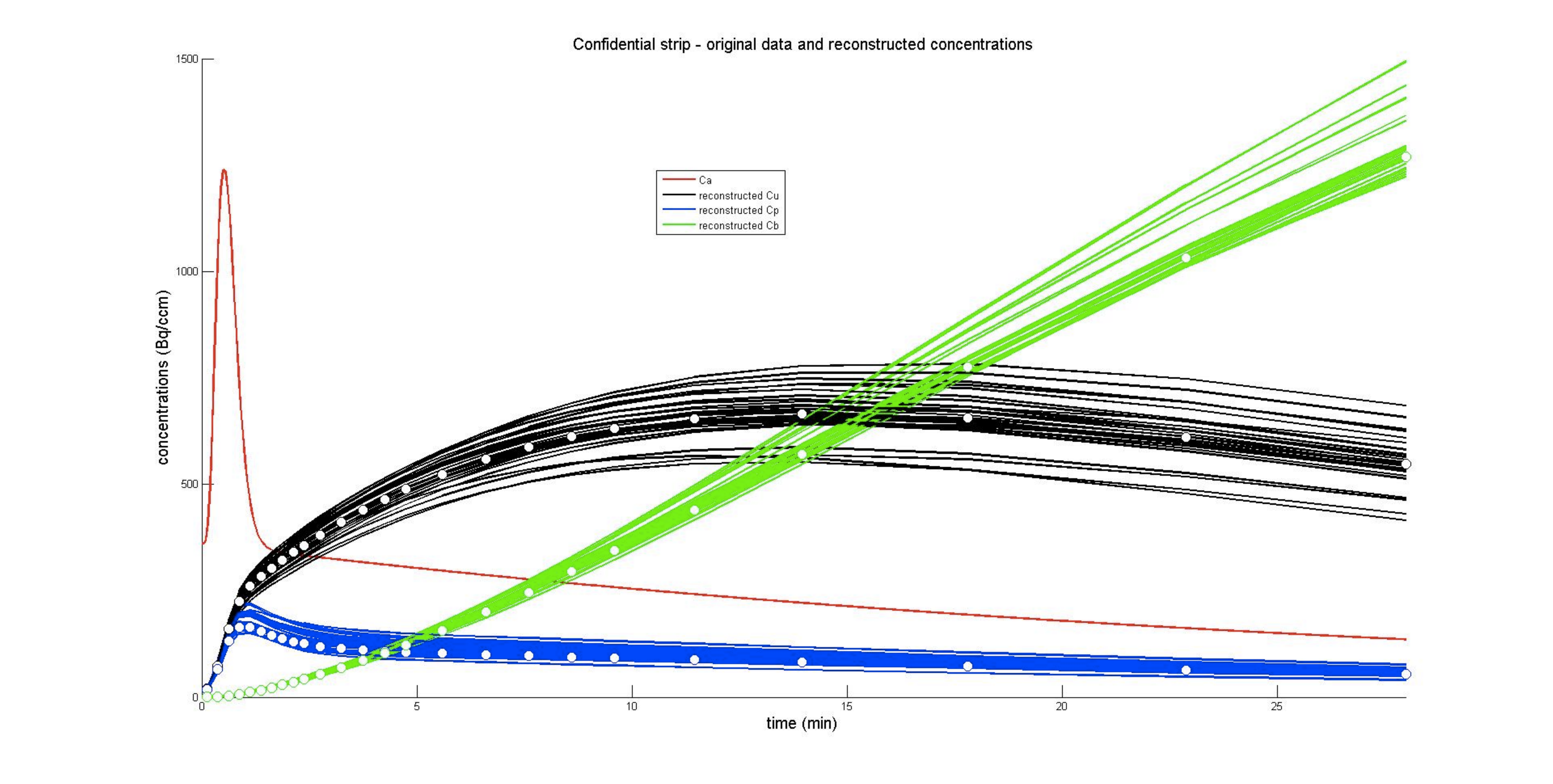,width=6.5cm} 
\end{tabular}
\caption{Results obtained with the following ACO parameters: $P=13$, $Q=7$, $q=0.015$, $\xi=0.4$ for 30 runs of the algorithm. Left panel: red line represents $C_b$, blue line represents $C_t$,  black line represents $C_p$; green line represents $C_u$. Right panel: superimposition of synthetic data (white dots) and reconstructed confidence strips of concentrations.}
\label{fig:simulated1}
\end{center}
\end{figure}

In order to produce the synthetic data we started from six initial values for the tracer coefficients. These selected values generate a matrix $\mathbf{A}$ and a vector ${\mathbf{K}}$ as in (\ref{system-4}) and (\ref{system-5}). The corresponding solution for $C_{t}$, $C_{p}$ and $C_{u}$ associated to the Cauchy problem (\ref{system-1}) are given in (\ref{solution-full-1}), (\ref{solution-full-2}) and (\ref{solution-full-6}). The solutions, represented in Figure \ref{fig:simulated1}, left panel, are sampled at 27 time points corresponding to the distribution of acquisition times experimentally performed by 'Albira'. The red line in the panel represents the TAC $C_{b}$ that has been obtained by fitting with a gamma variate function a set of real measurements acquired from a healthy mouse in a very controlled experiment \cite{Golish}. The vectors corresponding to the discretization of $C_{t}$ and $C_{p}$ are then summed together and affected by Poisson noise; the same kind of noise is applied to the vector corresponding to the discretization of $C_{u}$. ACO is applied 30 times to these synthetic data set to obtain 30 sets of reconstructed values of the exchange coefficients; in correspondence with each set we then computed the concentration curves and super-imposed them in Figure \ref{fig:simulated1}, right panel, thus producing confidence strips for the reconstructed concentrations. In Table \ref{table:results1} we present the average values of the tracer coefficients and the standard deviations of the values reconstructed by the 30 runs of ACO over the same set of input vectors.

\begin{table}
\begin{center}
\begin{tabular}{c|c|c|c|c|c|c}
\hline
&$k_{bt}$ & $k_{tp}$ & $k_{pt}$ & $k_{up}$ & $k_{tb}$ & $k_{pb}$ \\
\hline 
ground truth & 1  & 0.02  &  0.02 &   0.08 &   0.3&    0.3 \\
\hline
Mean  &1.0164   & 0.0201  &  0.0219 &   0.0819 &   0.3306&    0.3040 \\
\hline 
Std &    0.1060  &  0.0094 &   0.0075  &  0.0049  & 0.0341   & 0.0199\\
\hline 
\end{tabular}
\end{center}
\caption{Simulated values of tracer coefficients; average values and standard deviations over 30 runs of ACO.}
\label{table:results1}
\end{table}

\begin{figure}
\begin{center}
\begin{tabular}{cc}
\psfig{figure=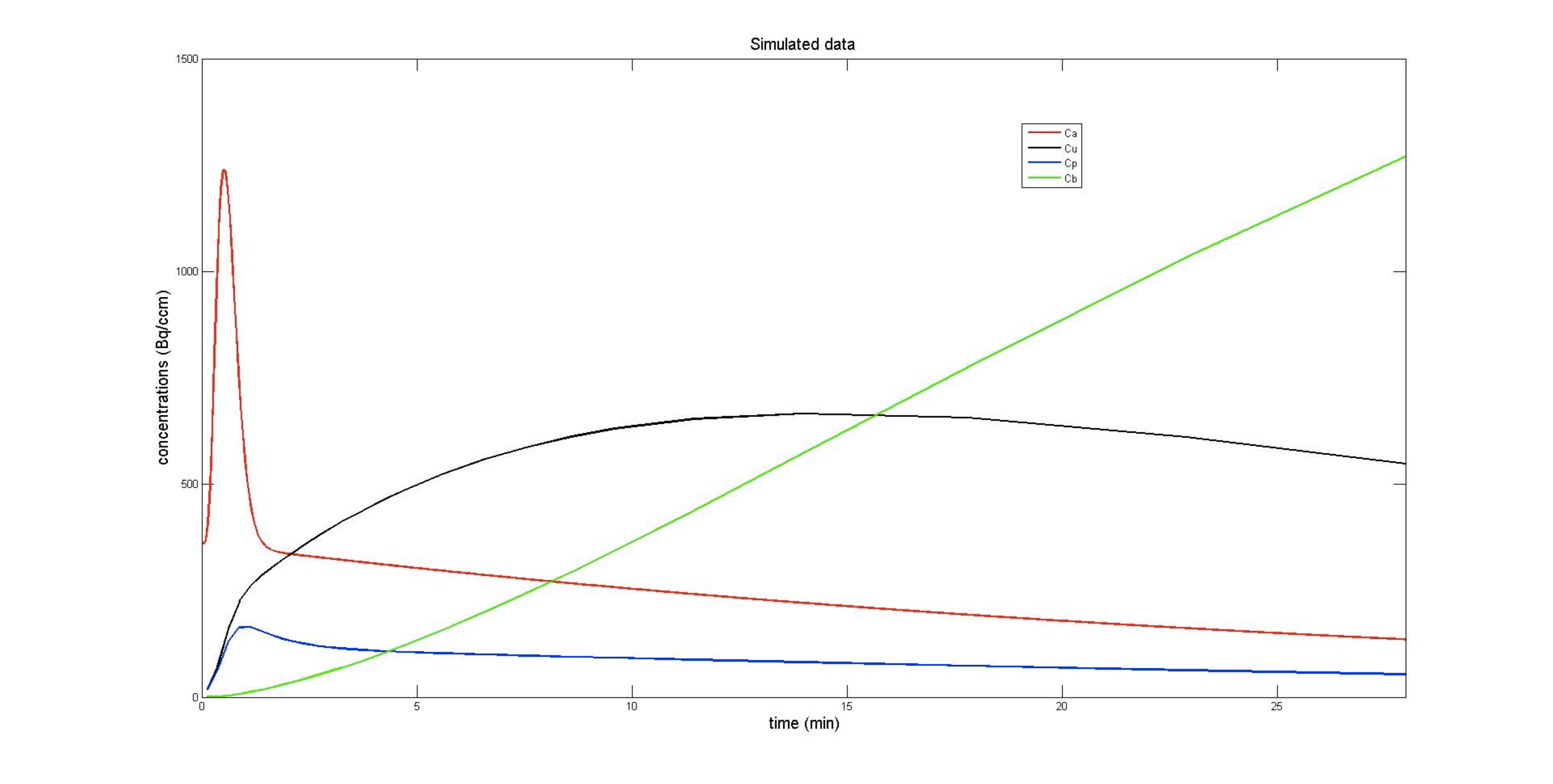,width=6.5cm} &
\psfig{figure=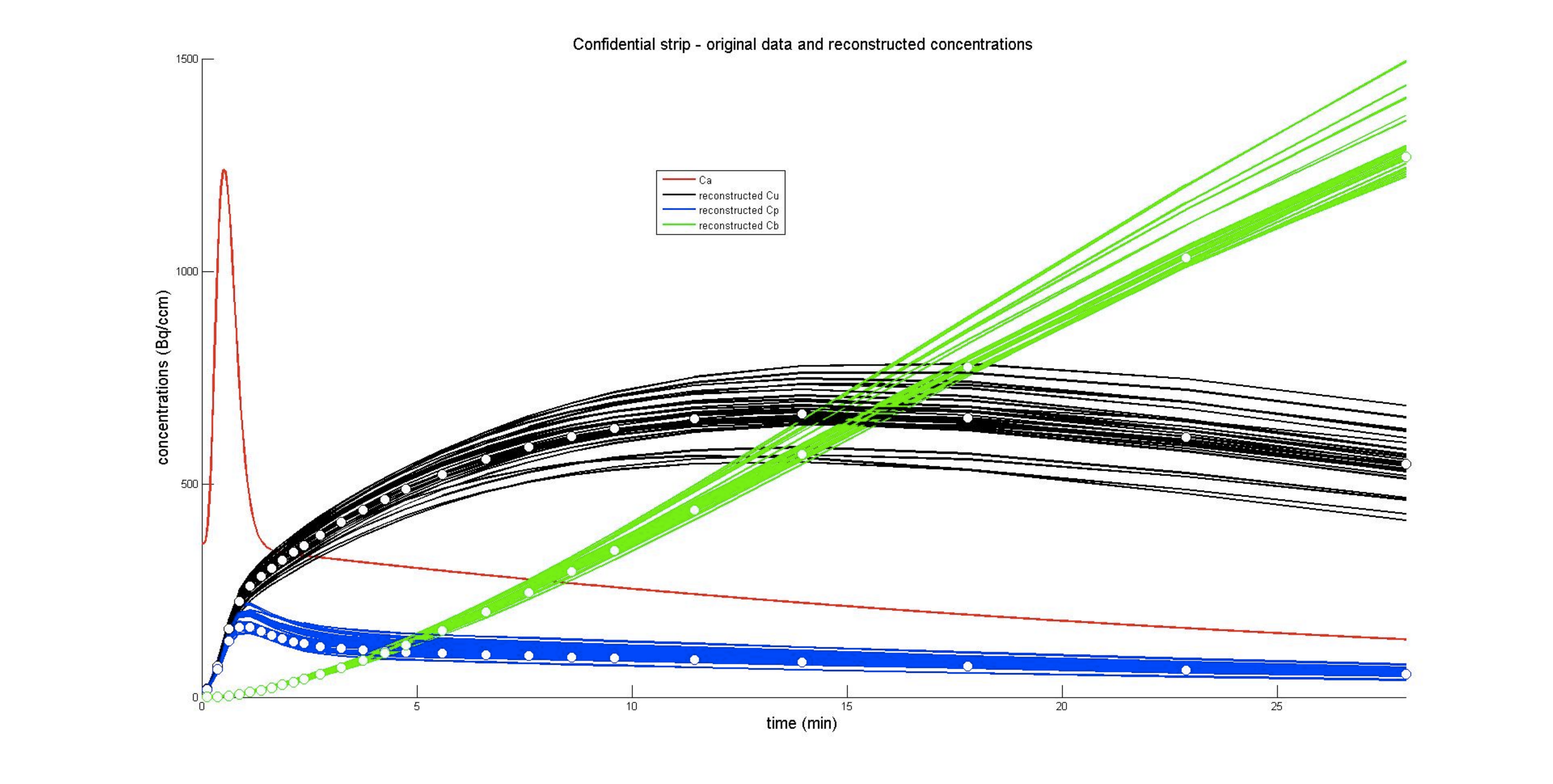,width=6.5cm} 
\end{tabular}
\caption{Results obtained with the following ACO parameters: $P=13$, $Q=7$, $q=0.015$, $\xi=0.4$ for 30 runs of the algorithm. Left panel: red line represents $C_b$, blue line represents $C_t$,  black line represent $C_p$ while green line represent $C_u$. Right panel: superimposition of synthetic data (white dots) and reconstructed confidence strips of concentrations.}
\label{fig:simulated2}
\end{center}
\end{figure}

The same kind of experiment has been performed again, this time choosing as initial values the ones that give rise to an upper triangular system matrix as in (\ref{upper-triangular}). The results for this case are in Figure \ref{fig:simulated2} and Table \ref{table:results2}. The reason why, in this case, ACO gives zero to the average and standard deviation of $k_{tp}$ is due to the fact that our implementation contains a threshold for the parameter outputs (equal to $10^{-3}$). In this experiment $k_{tp}$ is under-threshold for all runs and therefore the average and standard deviation are set to null.

We agree that in these tests the procedure for generating the synthetic concentrations and the one for reconstructing the tracer coefficients from them are based on the same equations (in a sort of 'inverse crime' procedure). However the synthetic data are affected by Poisson noise and in any case the aim of these two numerical applications was simply to validate the reliability and stability of ACO when applied, for the first time, to a compartmental analysis problem.

\begin{table}
\begin{center}
\begin{tabular}{c|c|c|c|c|c|c}
\hline
&$k_{bt}$ & $k_{tp}$ & $k_{pt}$ & $k_{up}$ & $k_{tb}$ & $k_{pb}$ \\
\hline 
ground truth & 0.8  & 0  &  0.02 &   0.08 &   0.4&    0.2 \\
\hline
Mean  &    0.8921    &     0   & 0.0209  &  0.1089    &0.4121    &0.1979\\
\hline
Std&    0.1078       &  0    &0.0044  &  0.0102 &   0.0526  &  0.0166\\
\hline 
\end{tabular}
\end{center}
\caption{Simulated values of tracer coefficients; average values and standard deviations over 30 runs of ACO.}
\label{table:results2}
\end{table}

\begin{table}
\begin{center}
\begin{tabular}{c|c|c|c|c|c|c}
\hline
&$k_{bt}$ & $k_{tp}$ & $k_{pt}$ & $k_{up}$ & $k_{tb}$ & $k_{pb}$ \\
\hline 
initial values &0.9572  &  0.4854 &   0.8003    &0.1419  &  0.4218  & 0.9157 \\
\hline 
mean  &  1.10  &  0.03   & 0.04   & 0.30  &  0.26  &  0.25\\
\hline
std&   0.41 &   0.02  &  0.03   & 0.06   & 0.02   & 0.02\\
\hline
\end{tabular}
\end{center}
\caption{Initial random values of tracer coefficients; average values and standard deviations over 20 runs of ACO.}
\label{table:results3}
\end{table}

\subsection{Application to real measurements}
\begin{figure}
\begin{center}
\begin{tabular}{cc}
\psfig{figure=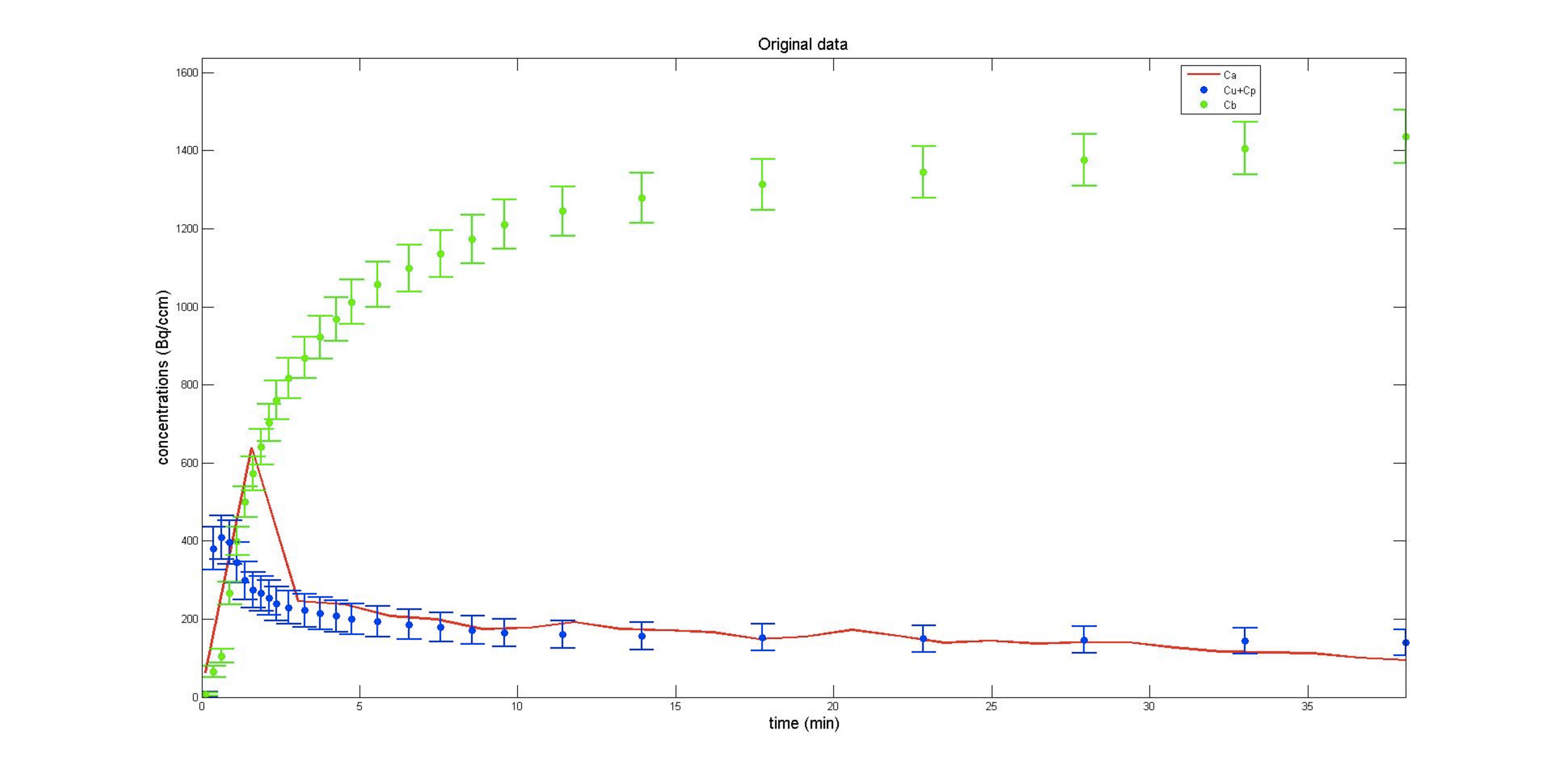,width=6.5cm} &
\psfig{figure=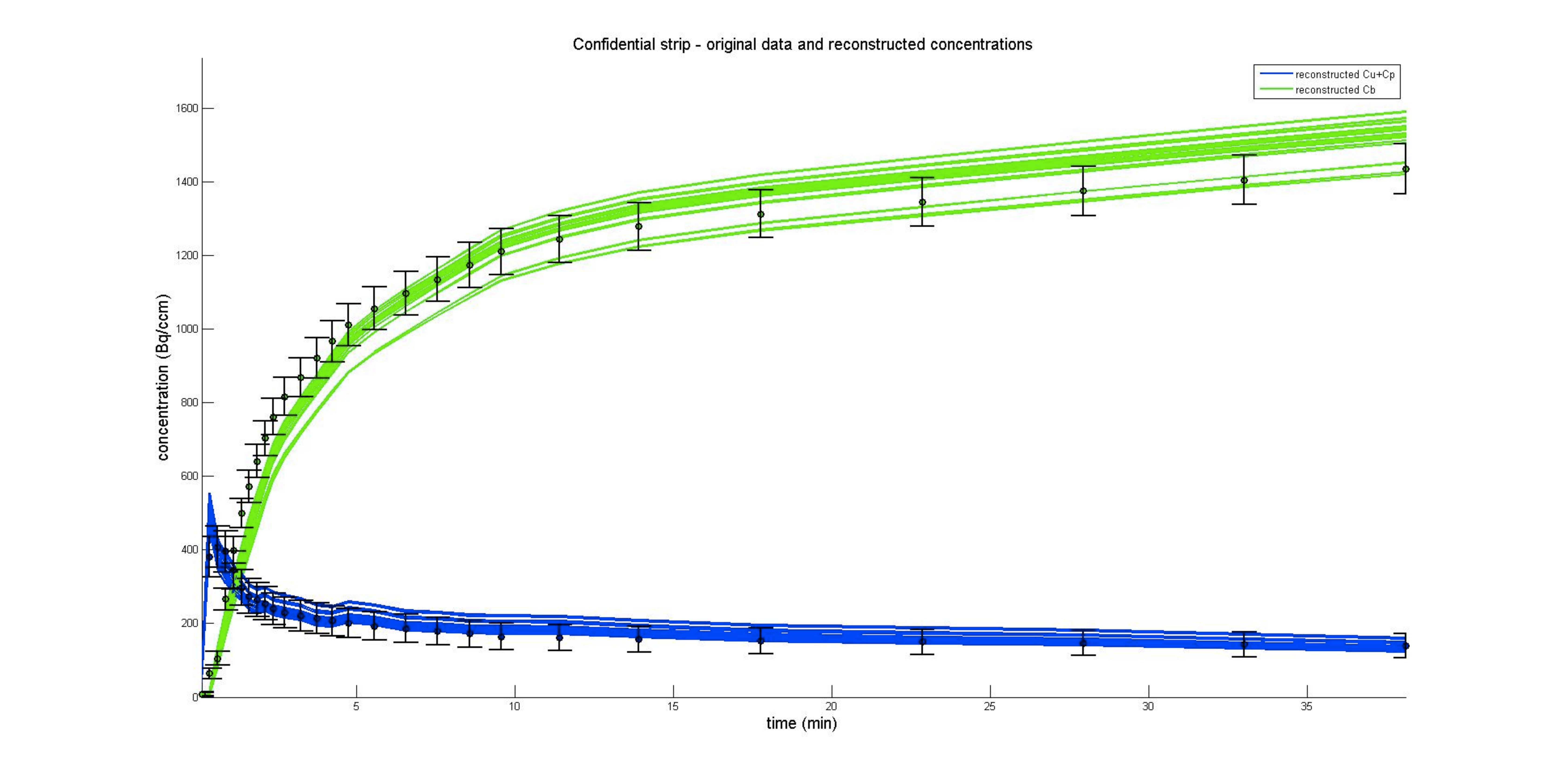,width=6.5cm} 
\end{tabular}
\caption{Results obtained with the following ACO parameters: $P=25$, $Q=13$, $q=0.0001$, $\xi=0.65$ for 20 runs of ACO. Left panel: red line represents $C_{b}$, blue dots $C_{t}+C_{p}$ while green dots represent $C_{u}$. Right panel: superimposition of concentrations in the bladder (green) and in the kidneys (blue) computed by solving the forward problem where the tracer coefficients are reconstructed by ACO.}
\label{fig:real1}
\end{center}
\end{figure}

\begin{figure}
\begin{center}
\begin{tabular}{cc}
\psfig{figure=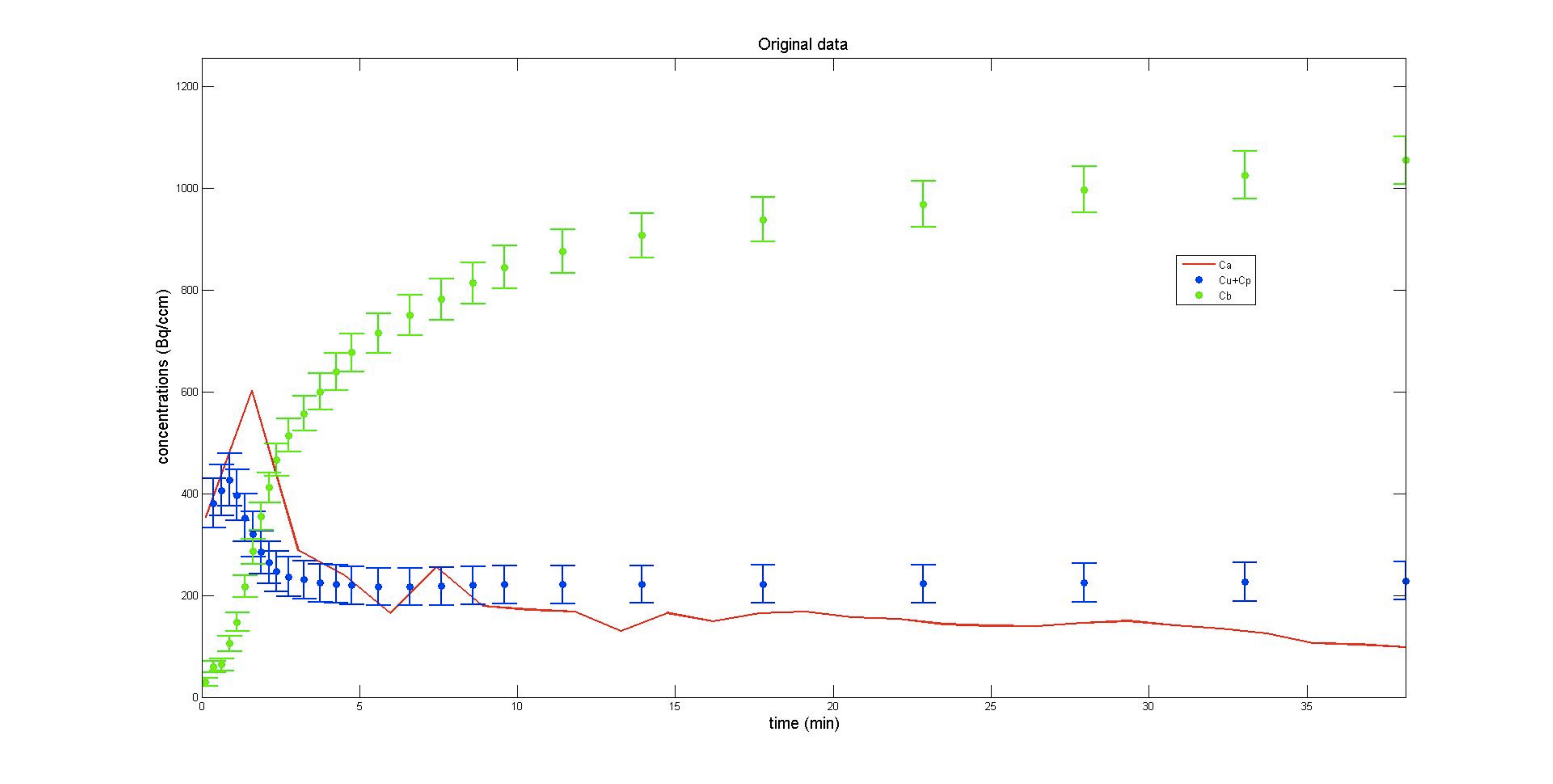,width=6.5cm} &
\psfig{figure=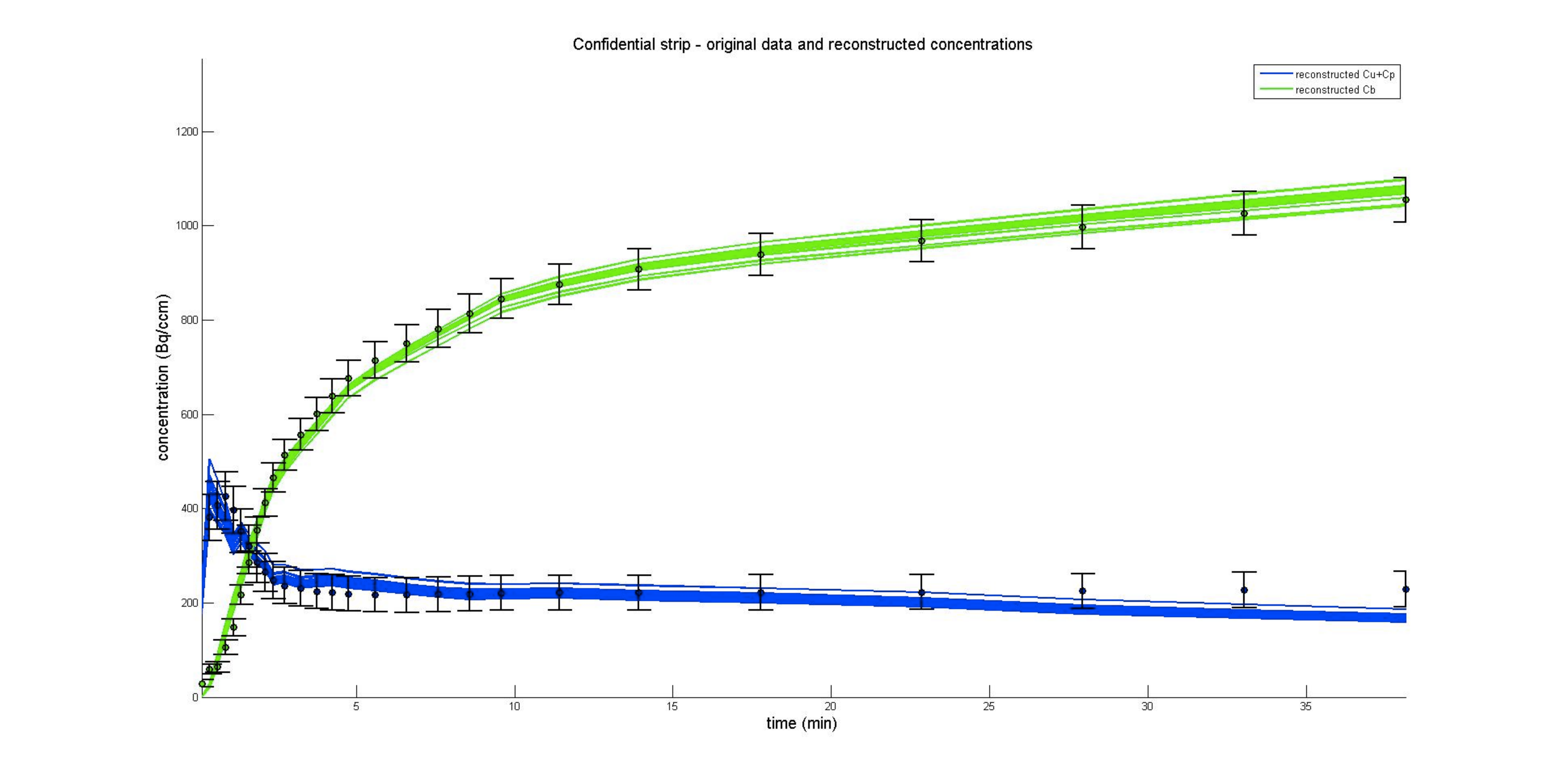,width=6.5cm} 
\end{tabular}
\caption{Results obtained with the following ACO parameters: $P=25$, $Q=13$, $q=0.001$, $\xi=0.45$ for 20 runs of ACO. Left panel: red line represents $C_{b}$, blue dots $C_{t}+C_{p}$ while green dots represent $C_{u}$. Right panel: superimposition of concentrations in the bladder (green) and in the kidneys (blue) computed by solving the forward problem where the tracer coefficients are reconstructed by ACO.}
\label{fig:real2}
\end{center}
\end{figure}

For the first experiment we considered a healthy murine model injected with FDG and acquired the corresponding activity by means of a dynamic acquisition paradigm over 27 experimental time points. The images have been reconstructed by applying an Expectation-Maximization iterative algorithm \cite{Dempster}
and ROIs have been drawn on the reconstructed images around the left ventricle to reproduce the time activity curve. ROIs have been also drawn around the kidneys and the bladder in order to compute the input concentrations. In Figure \ref{fig:real1}, left panel, the red line describes the TAC $C_{b}$ (we have plotted the solid line connecting the measured concentrations in order to distinguish the input function from the other concentrations). The green points represent the concentrations for the bladder. In order to compute the error bars for the measurements, the concentrations have been multiplied times the volume of the bladder and the resulting radioactivity is characterized by Poisson statistics, so that an estimate of the standard deviation is provided by the square root of the activity. The error bars of the concentration are therefore obtained by normalizing this square root with volume. The blue points correspond to the measured concentration in the kidneys and the error bars are computed as before. Then ACO has been applied 20 times against these data. The initialization values for the tracer coefficients are the same for all 20 runs and are obtained by drawing randomly a 6-ple within $(0,1)^6$. The results of this analysis are given in Table \ref{table:results3} containing the average and standard deviation (over the 20 realizations) of the tracer coefficients. Then, for each run, the set of reconstructed tracer coefficients are used to solve the direct problem in order to obtain reconstructions of $C_{u}$ and $C_{p} + C_{t}$. The confidence strip resulting from the superposition of the 20 reconstructions of the concentration are represented in Figure \ref{fig:real1}, right panel. A similar analysis has been performed in the case of a murine model affected by breast cancer. The results of this second case are given in Figure \ref{fig:real2} and in Table \ref{table:results4}.

\begin{table}
\begin{center}
\begin{tabular}{c|c|c|c|c|c|c}
\hline
&$k_{bt}$ & $k_{tp}$ & $k_{pt}$ & $k_{up}$ & $k_{tb}$ & $k_{pb}$ \\
\hline 
initial values &0.1988  &  0.0916  &  0.6170  &  0.7276  &  0.9138  &  0.2503\\
\hline
mean &    0.95   & 0.07&    0.05  &  0.18  &  0.19 &   0.19\\
\hline
std & 0.17  &  0.05 &   0.02  &  0.01 &   0.01  &  0.02\\
\hline
\end{tabular}
\caption{Initial random values of tracer coefficients and their average values and standard deviations over 20 runs of ACO.}
\label{table:results4}
\end{center}
\end{table}

\section{Comments and conclusions}

This paper deals with renal flow of a radioactive tracer, [$^{18}$F]-FDG, injected into a mouse. The time evolution of tracer concentrations inside kidney and from kidney to bladder has been modeled by a linear system of ordinary differential equations with constant coefficients. The time variation of the total concentration of activity inside kidney and bladder (essentially, the sum of the solutions) has been estimated through an analysis of micro-PET data. The six constant exchange coefficients, which provide information on system's physiology and metabolism, have been regarded as unknowns. The related inverse problem has been solved by applying a statistical optimization procedure, referred to by the acronym ACO. Resulting applications to real and synthetic data have been shown and discussed.    

The mathematical approach described in this paper provides estimates of the six unknown coefficients. Unlike techniques based on graphical analysis, it does not require any distinction between irreversible or reversible uptake of tracer, nor identification of a time value after which suitable expressions evaluated from the data become linear in time \cite{Schmidt}. Moreover, the graphical methods provide fewer parameters, usually slopes and intercepts, which can be interpreted as functions of the original model parameters \cite{Logan}.  Technically, the general character of the optimization procedure based on ACO makes it applicable to compartmental model structures of a variety of types, provided measured data on the related total concentrations of activity are available. For example, no practical restriction on the number of compartments involved is required if the direct problems involved in the ACO procedure are solved numerically, instead of finding the analytic representation of the solutions, as we have done in the present paper. Similarly, the fully numeric approach does not require limitations such as the constraint that the system matrix is sign-symmetric \cite{Hearon,Schmidt-2} (indeed, the upper and lower triangular matrices considered in Section 2 do not comply with sign-simmetry). We have also seen that application of ACO to synthetic data corresponding to an upper triangular matrix leads to the final estimate $k_{tp}=0$. This shows that the ACO approach is also capable of recovering vanishing rate constants. 

The physiological basis for this study relies on the broad utilization of FDG in the diagnosing and staging of cancer. In fasting patients, this tracer accurately maps the insulin independent glucose metabolism as an index of aggressiveness and growth rate of neoplastic lesions. From a theoretical point of view, dynamic PET imaging might provide a quantitative estimation of glucose consumption in cancer by measuring tracer rate constants and serum glucose levels. However, the long acquisition time requested by this approach conflicts with the routine activity since it cannot be applied to the whole body, it is hardly tolerated by patients and it markedly reduces the number of diagnostic procedures. To at least partially overcome this limitation, current clinical protocols imply that the whole body is imaged in a single shot late after tracer injection. Under this condition, glucose consumption is approximated by the measurement of tracer accumulation, normalizing for administered radioactivity and distribution mass. However, although this standardized uptake value (SUV) is usually considered an adequate surrogate of cancer glucose uptake \cite{Wahl} it also and equally depends upon tracer bioavailability. Most often this latter variable is neglected since blood FDG clearance by tissues is assumed to be relatively stable in the different patients.
Nevertheless, differently from glucose, a significant amount of FDG is excreted in the urines. This tracer sequestration has been already documented and attributed to the low affinity of glucose symporters dedicated to glucose reabsorption from preurine in the kidney tubules \cite{Moran,Qiao}. Our microPET approach offered the unique advantage to simultaneously measure the time concentration curves defining arterial input function (measured in the left ventricle), kidney tracer handling and urinary sequestration in the bladder without any need to evaluate the complex FDG transit in renal pelvis \cite{Kosuda} and its uncertainties. This direct estimation documented that tracer exchange from preurine to kidney parenchima is extremely low. Kidney glucose reabsorption is an active process by two Na+-dependent symporters present in the apical membrane of tubular cells (SGLT-1 and SGLT-2). The hydroxyl groups in D-glucose play a critical role in the recognition and transportation of glucose out of the proximal tubules of the kidney. Removal of the hydroxyl group at the 2 position of D-glucose diminishes the affinity of this molecule for both the forms of SGLT \cite{Kleinzeller}.
Urinary excretion has been most often considered as a limitation for FDG in studying cancer of the urinary tract \cite{Jana}. However, the evidence of an extremely low reabsorption of tracer filtered in the preurine has a more general relevance: it implies that tracer loss caused by glomerular filtration markedly affects FDG persistence in the blood. Accordingly, SUV might largely differ in cancer lesions characterized by a similar metabolism being lower in patients with preserved glomerular filtration rate and higher in the presence of an impaired renal function. The disposal of methods able to estimate the relationship between the easily available glomerular filtration rate and FDG bioavailability might thus improve the clinical standardization of tracer uptake as a marker for tumor aggressiveness, beyond the current standard application of SUV.

\end{document}